\documentclass[12pt]{amsart}
\usepackage{amsmath,amssymb,latexsym}
\usepackage{fullpage}
\usepackage{amscd}
\theoremstyle{plain}
\newtheorem{theorem}{Theorem}

\newtheorem{proposition}{Proposition}

\newtheorem{conjecture}{Conjecture}

\begin{document}

\title{Tensor Product $L$-Functions On Metaplectic Covering Groups of $GL_r$}
\author{David Ginzburg }
\address{Ginzburg: School of Mathematical Sciences, Tel Aviv University, Ramat Aviv, Tel Aviv 6997801,
Israel}
\email{ginzburg@post.tau.ac.il}

\thanks{The author is partly supported by the Israel Science
Foundation grant number  461/18.}

\begin{abstract}
In this note we compute some local unramified integrals defined on metaplectic  covering groups of $GL$. These local integrals which were introduced by Suzuki, represent the standard tensor product $L$ function $L(\pi^{(n)}\times \tau^{(n)},s)$ and extend the well known local integrals which represent $L(\pi\times \tau,s)$. The computation is done using a certain "generating function" which extends a similar  function introduced by the author in a previous paper. In the last section we discuss the Conjectures of Suzuki and introduce a global integral which unfolds to the above local integrals. This last part is mainly conjectural and relies heavily on the existence of Suzuki representations defined on covering groups. In the last subsection we introduce a new global doubling integral which represents the partial tensor product $L$ function $L^S(\pi\times \tau,s)$.

\end{abstract}
\maketitle
\section{Introduction}\label{intro}
Let $F$ denote a local nonarchimedean field which contains all $n$th roots of unity. Denote by $GL_r^{(n)}$ the $n$-fold metaplectic covering group of $GL_r$ defined over the field $F$. When $n=1$ we take the group $GL_r^{(n)}$ to be the linear group $GL_r$. 
Let $\pi^{(n)}$ denote an unramified representation of $GL_r^{(n)}$, and similarly, let $\tau^{(n)}$ denote an unramified representation of $GL_m^{(n)}$. When $n=1$, we denote by $\pi$ and $\tau$ the representations $\pi^{(n)}$ and $\tau^{(n)}$.

When $n=1$, and assuming that $r<m$, it follows from the work of \cite{J-PS-S} that the integral
\begin{equation}\label{ra1}
\int\limits_{V_r\backslash GL_r}W_\pi(g)\overline{W}_\tau\begin{pmatrix} g&\\ &I_{m-r}\end{pmatrix}|g|^{s-\frac{m-r}{2}}dg
\end{equation}
converges for $\text{Re}(s)$ large, and is equal to the local $L$ function $L(\pi\times\tau,s)$. When $n=1$ and $r=m$, the situation is slightly different. This case was also covered in \cite{J-PS-S}, and will be omitted here. In integral \eqref{ra1}, $W_\pi$ denotes the Whittaker function associated to the unramified vector in the space of the representation $\pi$, normalized to equal one at the identity. As is well known, if we assume $W_\pi(e)=1$, the function $W_\pi(g)$ is uniquely determined. A similar statement holds for $W_\tau$. Also, in integral \eqref{ra1}, the group $V_r$ is the maximal unipotent subgroup of $GL_r$.

Assume that $n>1$ and $r\le m$. The naive approach to extend integral \eqref{ra1} to
covering groups, by replacing the representations $\pi$ and $\tau$ by 
$\pi^{(n)}$ and $\tau^{(n)}$, does not work. We do not get the local $L$-function $L(\pi^{(n)}\times\tau^{(n)},s)$. Maybe the main reason for that is that in the covering groups the space of Whittaker functions attached to an unramified representation, has in many cases a dimension greater than one. 

The first to realize how to extend Rankin-Selberg constructions from linear groups to covering groups were  Bump and Hoffstein. See \cite{B-H1} and \cite{B-H2}. In these papers they studied some low rank cases. They realized that to make things work one has to consider convolutions of unramified representations with the Theta representation attached to the group $GL_n^{(n)}$. In the global and local setup, these representations were constructed for any  $n$ in \cite{K-P}, and they have the important property that they have a unique Whittaker function. Later, Bump and Friedberg, see \cite{B-F} studied the case when $r=1$. They 
wrote the local, and global, integral which they conjectured to equal the 
$L$-function $L(\tau^{(n)},s)$, and proved their conjecture for $m=2,3$ and any $n$. In \cite{G1} the Conjecture stated in \cite{B-F} was proved and hence we now have a global construction which represents the partial standard $L$ function $L^S(\tau^{(n)},s)$ attached to the cuspidal representation $\tau^{(n)}$.
Suzuki, in \cite{S}, studied the most general case. His paper contains a local integral, which can be viewed as an extension of integral \eqref{ra1}, and is conjectured to equal the local $L$-function $L(\pi^{(n)}\times\tau^{(n)},s)$. He proved it for some cases, stated in \cite{S} Section 9. These cases are when $r=1,2$ and $m<n$ if $r=1$, and $m<2n$ if $r=2$. His construction also involves a global integral which we will describe in the last section of this paper. We remark that it is clear that in both papers, \cite{B-F} and \cite{S}, the construction should work in complete generality. The problems were computational but not conceptual. 

Following \cite{S}, we now describe the general construction. Keeping the above notations, we assume that  $r< nm$. Following \cite{S} page 754, we consider the local integral
\begin{equation}\label{ra2}
\int\limits_{V_r\backslash GL_r}W_{\pi^{(n)}}(g)\overline{W}_{\tau^{(n)},nm}\begin{pmatrix} g&\\
&I_{nm-r}\end{pmatrix}|g|^{s-\frac{nm-r}{2}}dg
\end{equation}
The function $W_{\pi^{(n)}}(g)$ is {\sl{any}} Whittaker function associated with the unramified vector in the space of the representation $\pi^{(n)}$, normalized to equal one at the identity. As mentioned above, in many cases the dimension of the space of these functions is greater than one. To describe the function $W_{\tau^{(n)},nm}$, we first consider the induced representation
\begin{equation}\label{ind1}
\epsilon^{(n)}(\tau^{(n)})=Ind_{P_{n,m}^{(n)}}^{GL_{nm}^{(n)}}(\tau^{(n)}\otimes\tau^{(n)}\otimes\ldots\otimes\tau^{(n)})\delta_{P_{n,m}}^{\frac{nm-1}{2nm}}
\end{equation}
Here $P_{n,m}$ is the parabolic subgroup of $GL_{nm}$ whose Levi part is $GL_m\times\ldots\times GL_m$ where $GL_m$ appears $n$ times. It follows
from \cite{S} Section 3.3 that up to a constant, the unramified constituent of this induced representation has a unique Whittaker function. Let $W_{\tau^{(n)},nm}$ denote the Whittaker function associated with the unramified vector and normalized to be one at the identity. We emphasize that  the fact that the above unramified representation has a unique Whittaker function does not guarantee that the integral \eqref{ra2} is independent of the choice of the Whittaker function $W_{\pi^{(n)}}(g)$. 

Notice that integral \eqref{ra2} is  an extension of integral \eqref{ra1}. When $n=1$, the above induced representation is equal to $\tau$, and integral \eqref{ra2} equals to integral \eqref{ra1}.
Indeed, for $n=1$ we have $W_{\tau^{(n)},nm}=W_\tau$.

Our main Theorem is
\begin{theorem}\label{th1}
With the above notations, assume that  $r<nm$. Then
integral \eqref{ra2} is equal to the local $L$ function $L(\pi^{(n)}\times\tau^{(n)},ns-(n-1)/2)$.
\end{theorem}
For the definition of the local unramified tensor product $L$ function, see Section \ref{not} right after equation \eqref{ind2}.
The proof of this Theorem is a simple extension of the ideas used in \cite{G1}. The idea is to introduce a function, denoted by $W_{\tau^{(n)},nrm}^{(n)}(h)$ defined on the group $GL_{nrm}^{(n)}$.  We shall do that in Section \ref{not}. This function is not a Whittaker function, but still have nice left invariant properties. Let $K_r$ denote the maximal compact subgroup of $GL_r$. If we assume that $|n|_F=1$ then $K_r$ splits in $GL_r^{(n)}$.
Then we consider the integral
\begin{equation}\label{ra3}
I=\int\limits_{ GL_r}\omega_{\pi^{(n)}}(g)\overline{W}_{\tau^{(n)},nrm}^{(n)}\begin{pmatrix} g&\\
&I_{r(nm-1)}\end{pmatrix}|g|^{s'}dg
\end{equation}
where $s'=s-(mnr-2r+1)/2$ . Here $\omega_{\pi^{(n)}}$ is the unique spherical function defined on $GL_r^{(n)}$ and is attached to $\pi^{(n)}$. See \cite{B-F}. Thus, $\omega_{\pi^{(n)}}$ is bi-$K_r$ invariant. In particular we have that $W_{\tau^{(n)},nrm}^{(n)}(\text{diag}(g,I_{r(nm-1)}))$ is also bi-$K_r$ invariant. We mention that when $m=1$, the function  $W_{\tau^{(n)},nrm}^{(n)}$ is the function $W_{nr}^{(n)}$ introduced in \cite{G1} equation (7). 

To prove Theorem \ref{th1}, we compute integral \eqref{ra3} in two different ways. First, using the uniqueness of the spherical function, we obtain the identity $\int f_{\pi^{(n)}}(kg)dk=
\omega_{\pi^{(n)}}(g)$. Here $f_{\pi^{(n)}}$ is the unramified vector in the space of  $\pi^{(n)}$ normalized to equal one at the identity, and the integration is over $K_r$. Plugging this into integral \eqref{ra3}, we obtain 
\begin{equation}\label{ra4}
I=\int\limits_{ GL_r}f_{\pi^{(n)}}(g)\overline{W}_{\tau^{(n)},nrm}^{(n)}\begin{pmatrix} g&\\
&I_{r(nm-1)}\end{pmatrix}|g|^{s'}dg
\end{equation}
Computing integral \eqref{ra4} we show that it is equal to $L(\pi^{(n)}\times\tau^{(n)},ns-(n-1)/2)$. This we do in Section \ref{first}. We mention that a crucial step in this computation is the use of Theorem 46  in \cite{K} which is based on Theorem 8.1 in \cite{C2}. See Section \ref{first} right after equation \eqref{fir10}. A second important remark is that integral \eqref{ra3} and the computation we perform in Section \ref{first} hold for all numbers $r,n$ and $m$ such that $mn>1$. In other words, the assumption that $r<nm$ is not needed here.

Similarly,  using the identity $\int W_{\pi^{(n)}}(kg)dk=\omega_{\pi^{(n)}}(g)$, we obtain
\begin{equation}\label{ra5}
I=\int\limits_{ GL_r} W_{\pi^{(n)}}(g)\overline{W}_{\tau^{(n)},nrm}^{(n)}\begin{pmatrix} g&\\
&I_{r(nm-1)}\end{pmatrix}|g|^{s'}dg
\end{equation}
Here $W_{\pi^{(n)}}$ is {\sl any} unramified Whittaker function associated with $\pi^{(n)}$, normalized to equal one at the identity.
Computing this integral we prove that it is equal to integral \eqref{ra2}. This we prove in Section \ref{second}. The combination of these two computations is the proof of Theorem \ref{th1}.

In Section \ref{global1} we discuss the global integrals which unfold to the corresponding local integrals \eqref{ra2}. When $r<m$, then as in the linear case these integrals are of Hecke type. When $r=m$, the integral we consider is a variation of the global integral considered in \cite{B-F}. To construct these global integrals we need to know that  certain representations, denoted by $\epsilon^{(n)}(\tau)$, are not zero. See Subsection \ref{suz1} for the definition of these representations.
This is not known in general, which implies that most of the discussion in that Section is still conjectural. 

In Subsection \ref{suz1} we recall the Suzuki conjecture as stated in \cite{S}. This conjecture states that the representations $\epsilon^{(n)}(\tau)$ exists, and
state their properties. Roughly speaking the Suzuki Conjecture states that for every irreducible cuspidal representation $\tau$ of $GL_m({\bf A})$ one can attach an automorphic representation $\epsilon^{(n)}(\tau)$ defined on the group $GL_{nm}^{(n)}({\bf A})$. Depending on $\tau$, the representation $\epsilon^{(n)}(\tau)$ is either a cuspidal representation, or a residue of an Eisenstein series.
  
In Subsection \ref{suz2}, assuming the Suzuki Conjecture, we introduce the global integrals. Thus, our main result in that Section is
\begin{theorem}\label{mainth}
Let $\tau$ denote an irreducible cuspidal representation of $GL_m({\bf A})$. Assume that the Suzuki Conjecture, as stated in Conjecture \ref{conj1}, holds for $\tau$. Let $\pi^{(n)}$ denote an irreducible cuspidal representation of $GL_r^{(n)}({\bf A})$. If $r\ne m$, then the partial $L$ function $L^S(\pi^{(n)}\times \tau,s)$ is holomorphic for all $s$. If $r=m$, then $L^S(\pi^{(n)}\times \tau,s)$ is holomorphic for all $s$, except possibly a simple pole at $s=0,1$. If this partial $L$ function has a pole then $\tau$ is in the image of the Shimura lift of an irreducible cuspidal representation $\tau^{(n)}$ defined on $GL_m^{(n)}({\bf A})$, and $\pi^{(n)}$ is isomorphic to the representation  $\widehat{\tau}^{(n)}$.
\end{theorem}
Here, the representation $\widehat{\tau}^{(n)}$ is defined as follows. For all functions $\varphi_{{\tau}^{(n)}}(h)$ in the space of ${\tau}^{(n)}$, we denote by  $\widehat{\tau}^{(n)}$ the representation of $GL_m^{(n)}({\bf A})$ generated
by all functions $\varphi_{{\tau}^{(n)}}(^th^{-1})$. By the Shimura lift we mean the following. Given an irreducible cuspidal representation $\tau^{(n)}$ defined on the group $GL_m^{(n)}({\bf A})$, there is an irreducible cuspidal representation $\tau$ defined on $GL_m({\bf A})$ with the following property. The unramified parameters of every unramified constituent of the representation $\tau$  are obtained from the corresponding unramified parameters
of $\tau^{(n)}$ raised  to the power of $n$. This lift is proved only in some cases. For the group $GL_2$ it is was proved in \cite{F}. For the group $GL_3^{(2)}$
some initial result indicates that the conjecture can be proved using the minimal representation of $F_4^{(2)}({\bf A})$. See \cite{G2}.

Finally, in the last sub-section we introduce a new global doubling integral which represents the partial $L$ function $L^S(\pi\times \tau,s)$. This integral is similar in nature to the integral which unfolds to the Whittaker integral introduced in sub-section \ref{suz2} in equation \eqref{glob7}. The convergence of this integral is proved in a similar way as for integral \eqref{glob7}. Therefore we will only show that it unfolds to a factoriazable integral whose local unramified integral is integral \eqref{ra3}. This global doubling integral is constructed for the linear group only. It is possible that a certain refinement of this integral can be constructed also for metaplectic representations.

\section{Notations}\label{not}
In this Section we introduce notations and basic definitions of the representations and functions that we use. The main reference is \cite{G1}. See also  \cite{K} and \cite{C-F-G-K} for some basic notations and definitions. Since the proof of the basic Theorem is similar to the proof of the Theorems in \cite{G1}, we use here similar notations.

Let $\tau^{(n)}$ denote an unramified representation defined on the group $GL_m^{(n)}$.  Thus, we have $\tau^{(n)}=Ind_{B_m^{(n)}}^{GL_m^{(n)}}\mu\delta_B^{1/2}$. The precise definition of this representation is given in \cite{G1} Section 2. Here, we denote by $B_m$ the Borel subgroup of $GL_m$ consisting of all upper triangular matrices. We denote by $T_m$ the subgroup of $GL_m$ consisting of all diagonal matrices. The character $\mu$ is attached to the $m$ characters $\mu_1,\ldots,\mu_m$ of $F^*$, and we will assume that $\mu$ is in general position.  In a similar way we let $\pi^{(n)}$ denote the unramified representation $Ind_{B_r^{(n)}}^{GL_r^{(n)}}\chi\delta_B^{1/2}$. If the character $\chi$ is attached to the $r$ unramified characters $\chi_1,\ldots,\chi_r$, then the local unramified tensor $L$ function is defined by
\begin{equation}\label{ind0}
L(\pi^{(n)}\times\tau^{(n)},s)=\prod_{i,j}\frac{1}{(1-\chi_i^n(p)\mu_j^n(p)q^{-s})}
\end{equation}
Here $p$ is a generator of the maximal ideal in the ring of integers of $F$, and $q^{-1}=|p|_F$.

The function $W_{\tau^{(n)},nrm}^{(n)}$ is defined  in \cite{C-F-G-K} Section 2.1 and in \cite{K} Section 2.2. We will give here the definition at a local unramified place. 

First consider the unramified representation
\begin{equation}\label{ind2}
Ind_{P_{nr,m}^{(n)}}^{GL_{nrm}^{(n)}}(\tau^{(n)}\otimes\tau^{(n)}\otimes\ldots\otimes\tau^{(n)})\delta_{P_{nr,m}}^{\frac{nm-1}{2nm}}
\end{equation}
Here  $P_{nr,m}$ is the standard parabolic subgroup of $GL_{nrm}$ whose Levi part is $GL_m\times\ldots\times GL_m$ where $GL_m$ appears $nr$ times. We denote by $\Theta(\tau^{(n)})$ the unramified sub-representation  of the representation \eqref{ind2}. 
 
Using induction by stages the representation $\Theta(\tau^{(n)})$ is also a sub-representation of 
\begin{equation}\label{ind3}
Ind_{Q_{nm,r}^{(n)}}^{GL_{nrm}^{(n)}}(\epsilon^{(n)}(\tau^{(n)})\otimes\ldots\otimes \epsilon^{(n)}(\tau^{(n)}))\delta_{Q_{nm,r}}^{\frac{nm-1}{2nm}}
\end{equation}
where $\epsilon^{(n)}(\tau^{(n)})$ was defined in \eqref{ind1}. Here $Q_{nm,r}$ is the parabolic subgroup of $GL_{nrm}$ whose Levi part is $GL_{nm}\times\ldots\times GL_{nm}$ where $GL_{nm}$ appears $r$ times.

Let $P_{nm,r}$ denote the parabolic subgroup of $GL_{nrm}$ whose Levi part is $GL_r\times\ldots\times GL_r$ where $GL_r$ appears $nm$ times. 
Let $U_{nm,r}$ denote its unipotent radical. It consists of all matrices of the form 
\begin{equation}\label{mat1}
\begin{pmatrix} I&X_{1,2}&X_{1,3}&\dots&X_{1,nm}\\
&I&X_{2,3}&\dots&X_{2,nm}\\
&&I&\ddots&\vdots\\ &&&\ddots&X_{nm-1,nm}\\ &&&&I\end{pmatrix}
\end{equation}
Here $I$ is the $r\times r$ identity  matrix, and $X_{i,j}\in 
Mat_{r\times r}$. Let $\psi$ denote a nontrivial  unramified character 
of $F$. Define a character
$\psi_{U_{nm,r}}$ of $U_{nm,r}$ as follows. For $u\in U_{nm,r}$ as above, define $\psi_{U_{nm,r}}(u)=\psi(\text{tr}(X_{1,2}+X_{2,3}+\cdots +
X_{nm-1,nm}))$. The stabilizer of $\psi_{U_{nm,r}}$ inside $GL_r\times GL_r \times\ldots\times GL_r$ is the group $GL_r^\Delta$ embedded diagonally. The embedding of $GL_r^\Delta$ inside $GL_{nrm}$ is given
by $g\mapsto \text{diag}(g,g,\ldots,g)$. 

As in \cite{C-F-G-K} and \cite{K}, to construct the function $W_{\tau^{(n)},nrm}^{(n)}$ we first define the Whittaker--Speh--Shalika functional $l$ on the space of $\Theta(\tau^{(n)})$. By definition this is a 
functional which satisfies the transformation property $l(\Theta(\tau^{(n)})(u)v)=\psi_{U_{nm,r}}(u)l(v)$ for all $u\in U_{nm,r}$ and all vector $v$ in $\Theta(\tau^{(n)})$. Assuming such a functional $l$ exists, we then define 
$W_{\tau^{(n)},nrm,v}^{(n)}(h)=l(\Theta(\tau^{(n)})(h)v)$. To construct such a functional $l$, we proceed as follows. 

Let $U_{nm,r}^0$ denote the subgroup of $U_{nm,r}$
consisting of all matrices as in \eqref{mat1} such that $X_{i,j}\in Mat^0_{r\times r}$. Here $Mat_{r\times r}^0$ is the subgroup of $Mat_{r\times r}$ consisting of all matrices $X$ such that $X[l_1,l_2]=0$ for all $l_1<l_2$, where $X[l_1,l_2]$ denotes the $(l_1,l_2)$-th entry of $X$. Notice that by restriction, $\psi_{U_{nm,r}}$ is also a character of 
$U_{nm,r}^0$.

Let $w_J=\text{diag}(J_{nm},J_{nm},\ldots,J_{nm})\in GL_{nmr}$. Here, $J_{nm}$ is the permutation matrix of size $nm$ which has ones on the other diagonal. In a similar way as in \cite{G1}, right after equation (7), let $w_0$ denote the Weyl element of $GL_{nrm}$ whose $(a+bnm, (a-1)r+b+1)$ entry is one, and all other entries are zeros. Here $1\le a\le nm$ and $0\le b\le r-1$.  
As in \cite{G1} equation (7), we define
\begin{equation}\label{local1}
W_{\tau^{(n)},nrm,f}^{(n)}(h)=\int\limits_{U_{nm,r}^0}f(w_Jw_0uh)\psi_{U_{nm,r}}(u)du
\end{equation}
Here $f$ is any vector in the space of $\Theta(\tau^{(n)})$. From the computation of the integral \eqref{fir1} done in Section \ref{first} below, it follows that if $f$ is the unramified vector in $\Theta(\tau^{(n)})$, then integral \eqref{local1} is not zero. In particular the space of such functions is not zero. 

It follows from the above, and from \cite{C1} Theorem 1.2,
\begin{proposition}\label{prop1}
The space of functions $W_{\tau^{(n)},nrm,f}^{(n)}$ defined as above on the representation $\Theta(\tau^{(n)})$ is one dimensional. Moreover, let $k^\Delta\in K_r\subset GL_r^\Delta$. Then $W_{\tau^{(n)},nrm,f}^{(n)}(k^\Delta h)=W_{\tau^{(n)},nrm,f}^{(n)}(h)$.
\end{proposition}
See  also \cite{G1} Proposition 1.

Finally, let $f=f_{nrm}^{(n)}$ denote the unramified function in the representation $\Theta(\tau^{(n)})$. Write $W_{\tau^{(n)},nrm}^{(n)}(h)$ for the function $W_{\tau^{(n)},nrm,f}^{(n)}(h)$. The properties of this function, which will be used in the next Sections are as follows. First,
given $u\in  U_{nm,r}$, we have $W_{\tau^{(n)},nrm}^{(n)}(uh)=
\psi_{U_{nm,r}}(u)W_{\tau^{(n)},nrm}^{(n)}(h)$. Then, for $k^\Delta\in 
K_r$ and for all $k\in K_{nrm}$ we have $W_{\tau^{(n)},nrm}^{(n)}(k^\Delta hk)=W_{\tau^{(n)},nrm}^{(n)}(h)$. Notice that this implies that for $g\in GL_r$, 
$W_{\tau^{(n)},nrm}^{(n)}(\text{diag}(g,I_{r(nm-1)}))$ is bi-$K_r$ invariant. In particular integral $I$ as defined in equation \eqref{ra3} is well defined.

\section{ First Computation of integral $I$}\label{first}
In this Section we compute the integral $I$ as given in \eqref{ra4}. We assume that $nm>1$. Thus, we compute the integral   
\begin{equation}\label{fir1}
I=\int\limits_{ GL_r}f_{\pi^{(n)}}(g)\overline{W}_{nrm}^{(n)}\begin{pmatrix} g&\\
&I_{r(nm-1)}\end{pmatrix}|g|^{s'}dg
\end{equation}
where  $s'=s-(mnr-2r+1)/2$ and $f_{\pi^{(n)}}$ is the unramified vector in the space of  $\pi^{(n)}$ normalized to equal one at the identity.
Plugging the expression for $\overline{W}_{\tau^{(n)},nrm}^{(n)}$ and performing the Iwasawa decomposition for $GL_r$, integral \eqref{fir1} is equal to
\begin{equation}\label{fir2}
I=\int\limits_{T_r}f_{\pi^{(n)}}(t)\int\limits_{V_r}
\int\limits_{U_{nm,r}^0}
\bar{f}_{nrm}^{(n)}
(w_Jw_0uv_0t_0)\psi_{U_{nm,r}}(u)dudv_0|t|^{s'}\delta_{B_r}^{-1}(t)dt
\end{equation}
where we recall that ${f}_{nrm}^{(n)}$ is the unramified vector in the representation $\Theta(\tau^{(n)})$.
In general, given $g\in GL_r$, we denote  $g_0=\text{diag}(g,I_{r(nm-1)})\in GL_{rnm}$. 

Let $U^1_{nm,r}$ denote the unipotent subgroup of $U_{nm,r}^0$ consisting of all matrices as in \eqref{mat1} such that $X_{1,2}[i,j]=0$ for all $i\ne j$. Then the inner integral in \eqref{fir2} over $V_r$ and $U_{nm,r}^0$ is equal to 
\begin{equation}\label{fir3}
\int\limits_{U_{nm,r}^1}
\bar{f}_{nrm}^{(n)}(w_Jw_0ut_0)\psi_{U_{nm,r}}(u)du
\end{equation}
The proof of this statement is exactly as the proof in \cite{G1} that integral (12) is equal to integral (13). We omit it here. 

Write $w_0U_{nm,r}^1w_0^{-1}=U_{nm,r}^2U_{nm,r}^3$. Here $U_{nm,r}^3$
consists of all matrices given in \cite{G1} equation (14) such that $Y_{i,j}\in Mat_{nm\times nm}$ and satisfy the conditions $Y_{i,j}[l_1,l_2]=Y_{i,j}[1,2]=0$ for all $l_1\ge l_2$.

The group $U_{nm,r}^2$ is the unipotent subgroup of $GL_{nrm}$ consisting of all upper unipotent matrices of the form $u_2=\text{diag}(v_{nm,1},v_{nm,2},\ldots,v_{nm,r})$ where $v_{nm,i}\in V_{nm}$. We recall that $V_{nm}$ is the maximal unipotent subgroup of $GL_{nm}$ consisting of upper unipotent matrices. Thus, the group $U_{nm,r}^2$ can be identified with $r$ copies of $V_{nm}$. Let $\psi_{V_{nm}}$ denote the Whittaker character of $V_{nm}$. This character is defined as in \cite{G1} equation (15) where $n$ is replaced by $nm$. Using this character, we define the character $\psi_{U_{nm,r}^2}$ of $U_{nm,r}^2$ by $\psi_{U_{nm,r}^2}(u_2)=\psi_{V_{nm}}(v_{nm,1})
\psi_{V_{nm}}(v_{nm,2})\ldots \psi_{V_{nm}}(v_{nm,r})$.
From this we deduce that integral \eqref{fir3} is equal to
\begin{equation}\label{fir4}
\int\limits_{U_{nm,r}^3}
f_{nrm,W}^{(n)}(u_3w_0t_0w_0^{-1})du_3
\end{equation}
where 
\begin{equation}\label{fir5}\notag
f_{nrm,W}^{(n)}(h)=
\int\limits_{U_{nm,r}^2}
\bar{f}_{nrm}^{(n)}(w_Ju_2h)\psi_{U_{nm,r}^2}(u_2)du_2
\end{equation}
Write $t=\text{diag}(a_1,a_2,\ldots,a_r)\in T_r$ where $a_i\in F^*$.
As in \cite{G1} after equation (16), we have $w_0t_0w_0^{-1}=\text{diag} (A_1,A_2,\ldots,A_r)$ where
$A_i=\text{diag}(a_i,I_{nm-1})$. Conjugating the matrix $w_0t_0w_0^{-1}$ to the left in integral \eqref{fir4} we obtain the factor
\begin{equation}\label{fir6}
\alpha(t)=(|a_2||a_3|^2|a_4|^3\ldots |a_r|^{r-1})^{nm-2}\notag
\end{equation}
from the change of variables in $U_{nm,r}^3$. Then, arguing as in \cite{G1} right after equation (17), we obtain that 
\begin{equation}\label{fir7}
I=\int\limits_{T_r}f_{\pi^{(n)}}(t)
\alpha(t)\bar{f}_{nrm,W}^{(n)}
(w_0t_0w_0^{-1})\delta_{B_r}^{-1}(t)|t|^{s'}dt
\end{equation}
Since the vector ${f}_{nrm,W}^{(n)}$ is the unramified vector in the space of the representation $\Theta(\tau^{(n)})$ which is the unramified sub-representation of \eqref{ind3}, then we have
\begin{equation}\label{fir8}\notag
\bar{f}_{nrm,W}^{(n)}(w_0t_0w_0^{-1})=\prod_{i=1}^rW_{\tau^{(n)},nm}
\begin{pmatrix} a_i&\\ &I_{nm-1}\end{pmatrix}\delta_{P_{nm,r}}^{\frac{nm-1}{2nm}}(w_0t_0w_0^{-1})
\end{equation}
Here, $W_{\tau^{(n)},nm}$ is the unique nontrivial Whittaker function associated with the unramified vector of the representation 
$\epsilon^{(n)}(\tau^{(n)})$ normalized to equal one at the identity. See right after \eqref{ra1}. Clearly, the right hand side is zero unless $|a_i|\le 1$. It follows from \cite{K} Theorem 46 that the right hand side is also zero unless $a_i=b_i^n$ for some $b_i\in F^*$. 

Recall
that $\pi^{(n)}=Ind_{B_r^{(n)}}^{GL_r^{(n)}}\chi\delta_{B_r}^{1/2}$. Since $t=\text{diag}(b_1^n,\ldots,b_r^n)$, we clearly have 
$f_{\pi^{(n)}}(t)=\chi(t)\delta_{B_r}^{1/2}(t)$. Hence, 
\begin{equation}\label{fir9}\notag
I=\int\limits_{T_r}\chi(t)
\alpha(t)\prod_{i=1}^rW_{\tau^{(n)},nm}
\begin{pmatrix} b_i^n&\\ &I_{nm-1}\end{pmatrix}
\delta_{P_{nm,r}}^{\frac{nm-1}{2nm}}(w_0t_0w_0^{-1})\delta_{B_r}^{-1/2}(t)|t|^{s'}dt
\end{equation}
We have 
$$\alpha(t)\delta_{P_{nm,r}}^{\frac{nm-1}{2nm}}(w_0t_0w_0^{-1})\delta_{B_r}^{-1/2}(t)=|b_1b_2\ldots b_r|^{n(nm-2)(r-1)/2}$$
Since $s'=s-(mnr-2r+1)/2$, we deduce that
\begin{equation}\label{fir10}
I=\prod_{i=1}^r\int\limits_{|b_i|\le 1}\chi_i(b_i^n)
W_{\tau^{(n)},nm}
\begin{pmatrix} b_i^n&\\ &I_{nm-1}\end{pmatrix}
|b_i^n|^{s-(nm-1)/2}db_i
\end{equation}
Finally, using Theorem 46 in \cite{K} and Theorem 8.1 in \cite{C2}, we obtain
\begin{equation}\label{fir11}\notag
I=\prod_{i=1}^rL(\chi_i^n\times\tau^{(n)},ns-(n-1)/2)=L(\pi^{(n)}\times\tau^{(n)} ,ns-(n-1)/2)
\end{equation}

\section{ Second Computation of integral $I$}\label{second}
In this Section we compute the integral $I$, using identity \eqref{ra5}. We assume that $r<nm$.
In other words we compute the integral
\begin{equation}\label{sec1}\notag
I=\int\limits_{ GL_r} W_{\pi^{(n)}}(g)\overline{W}_{\tau^{(n)},nrm}^{(n)}\begin{pmatrix} g&\\
&I_{r(nm-1)}\end{pmatrix}|g|^{s'}dg
\end{equation}
As before $s'=s-(mnr-2r+1)/2$. To do that we argue as in \cite{G1} Theorem 2. Factoring the measure over $GL_r$, we obtain that
\begin{equation}\label{sec2}
I=\int\limits_{ V_r\backslash GL_r} W_{\pi^{(n)}}(g)
\int\limits_{V_r}
\overline{W}_{\tau^{(n)},nrm}^{(n)}\begin{pmatrix} vg&\\
&I_{r(nm-1)}\end{pmatrix}\psi_{V_r}^{-1}(v)dv|g|^{s'}dg
\end{equation}

At this point we apply a similar proof as the proof of Theorem 2 in \cite{G1}. The only different in the proof is replacing $n$ there with $nm$. Thus, for all $g\in GL_r$, we obtain
 \begin{equation}\label{sec3}\notag
\int\limits_{V_r}
W_{\tau^{(n)},nrm}^{(n)}\begin{pmatrix} vg&\\
&I_{r(nm-1)}\end{pmatrix}\psi_{V_r}^{-1}(v)dv=
W_{\tau^{(n)},nm}^{(n)}
\begin{pmatrix} g&\\ &I_{nm-r}\end{pmatrix}
|g|^{(nm-1)(r-1)/2}
\end{equation}
Plugging this into equation \eqref{sec2}, we obtain the integral \eqref{ra2}.

\section{Global Constructions and the Conjecture of Suzuki}\label{global1}
\subsection{On the Conjecture of Suzuki}\label{suz1}
Let $\tau$ denote an irreducible cuspidal representation of $GL_m({\bf A})$.  The Conjecture of Suzuki, stated in the introduction of \cite{S}, attaches to $\tau$ an automorphic representation $\epsilon^{(n)}(\tau)$ of $GL_{nm}^{(n)}({\bf A})$. This is done as follows. 

First assume that $\tau$ is in the image of the conjectural Shimura lift of  an irreducible cuspidal representation $\tau^{(n)}$ of $GL_m^{(n)}({\bf A})$. Then, as in \cite{S} Section 8, one can form the multi-variable Eisenstein series $E_{\tau^{(n)}}^{(n)}(g,s_1,\ldots,s_n)$
defined on $GL_{nm}^{(n)}({\bf A})$ which is attached to the induced representation 
\begin{equation}\label{glob1}
Ind_{P_{n,m}^{(n)}({\bf A})}^{GL_{nm}^{(n)}({\bf A})}(\tau^{(n)}|\cdot|^{s_1}\otimes\tau^{(n)}|\cdot|^{s_2}\otimes\ldots\otimes\tau^{(n)}|\cdot|^{s_n})\delta_{P_{n,m}}^{1/2}
\end{equation}
Here $P_{n,m}$ is the parabolic subgroup of $GL_{nm}$ whose Levi part is $GL_m\times\ldots\times GL_m$ where $GL_m$ appears $n$ times. As shown in \cite{S} Section 8, this Eisenstein series has a simple pole at the
point $\sum_{i=1}^ns_i=0$ and $n(s_i-s_{i+1})=1$ for $1\le i\le n-1$. Let 
$\epsilon^{(n)}(\tau)$ denote the residue representation of 
$E_{\tau^{(n)}}^{(n)}(g,s_1,\ldots,s_n)$ at that point. 

We remark that one expects that for the existence of the residue representation it is not necessary to assume that $\tau^{(n)}$ has a Shimura lift to $\tau$. If, for example, we know that the partial $L$ function $L^S(\tau^{(n)}\times 
\widehat{\tau}^{(n)},s)$ has a simple pole at $s=1$, then the existence of the residue follows as in \cite{S}.  Unfortunately, at this stage it is not known how to establish in general the existence of the residue representation.

Write $\epsilon^{(n)}(\tau)=\otimes_\nu \epsilon^{(n)}(\tau)_\nu$. It is not hard to check that at an unramified non-archimedean place $\nu$, the local unramified constituent $\epsilon^{(n)}(\tau)_\nu$ is a sub-representation of
\begin{equation}\label{glob2}
Ind_{P_{n,m}^{(n)}(F_\nu)}^{GL_{nm}^{(n)}(F_\nu)}(\tau_\nu^{(n)}\otimes\tau_\nu^{(n)}\otimes\ldots\otimes\tau_\nu^{(n)})\delta_{P_{n,m}}^{\frac{nm-1}{2nm}}
\end{equation}
This local representation  was denoted by $\epsilon^{(n)}(\tau^{(n)})$ in equation \eqref{ind1}. 

Assume that the unramified parameters attached to $\tau_\nu$ are $\eta_{1,\nu},\ldots, \eta_{m,\nu}$. Then the parameters of $\tau^{(n)}_\nu$  are $\mu_{1,\nu},\ldots, \mu_{m,\nu}$ where $\mu_{i,\nu}^n=\eta_{i,\nu}$ for all $1\le i\le m$. Recall that unramified representations of $GL_m^{(n)}$ are determined uniquely by unramified characters on $n$ powers. Then, as explained in \cite{K}  section 2.2, we obtain that the unramified constituent $\epsilon^{(n)}(\tau)_\nu$ is a sub-representation of
\begin{equation}\label{glob3}
Ind_{P_{m,n}^{(n)}(F_\nu)}^{GL_{nm}^{(n)}(F_\nu)}(\Theta_{\mu_{1,\nu}}^{(n)}\otimes\Theta_{\mu_{2,\nu}}^{(n)}\otimes\ldots\otimes\Theta_{\mu_{m,\nu}}^{(n)})\delta_{P_{m,n}}^{\frac{1}{2}}
\end{equation}
Here $\Theta_{\mu_{i,\nu}}^{(n)}$ is the local Theta representation
of the group $GL_n^{(n)}(F_\nu)$ twisted by the character $\mu_{i,\nu}$. This representation was constructed in \cite{K-P}, and as proved there it has a unique Whittaker function normalized to equal one at the identity. Notice that in equation \eqref{glob2}, the parabolic subgroups we use for the induction process is $P_{n,m}$ whereas in \eqref{glob3} we used the group $P_{m,n}$. In general $P_{k,l}$ is the parabolic subgroup of 
$GL_{kl}$ whose Levi part is $GL_l\times\ldots\times GL_l$ where $GL_l$ appears $k$ times.

The second case to consider is when $\tau$ is not in the image of the Shimura lift from $GL_m^{(n)}({\bf A})$. In this case the Conjecture of Suzuki states that there is an irreducible cuspidal representation, denoted by $\epsilon^{(n)}(\tau)$, defined  on $GL_{nm}^{(n)}({\bf A})$,
and such that for almost all places its unramified constituent is the unramified constituent of the representation \eqref{glob3}. At each place $\nu$, the parameters $\mu_{i,\nu}$ are such that $\mu_{i,\nu}^n$ are the parameters of $\tau_\nu$.

We summarize this,
\begin{conjecture}\label{conj1}
Let $\tau$ denote an irreducible cuspidal representation of $GL_m({\bf A})$. Then to $\tau$ one can associate a representation $\epsilon^{(n)}(\tau)$, defined  on $GL_{nm}^{(n)}({\bf A})$, such that for almost all places, the unramified constituent of $\epsilon^{(n)}(\tau)$ at a place $\nu$ is the unramified constituent of the representation \eqref{glob3}. The representation 
$\epsilon^{(n)}(\tau)$ is generic, and its unramified constituent has a unique Whittaker function.
\end{conjecture}
This Conjecture can be described in the following diagram,
\begin{equation}\label{diag1}
\begin{matrix} &&(GL_{nm}^{(n)}({\bf A}),\epsilon^{(n)}(\tau))\\
&\nearrow&\\ (GL_m({\bf A}),\tau)&&\\ &\searrow&\\ &&( GL_{m}^{(n)}({\bf A}),\tau^{(n)})\end{matrix}
\end{equation}

When $n=1$ the Conjecture is trivial. When $n>1$, then parts of this Conjecture are proved. Clearly, if $\epsilon^{(n)}(\tau)$
is cuspidal then it is generic. If $\epsilon^{(n)}(\tau)$ is the residue representation then it follows from the Theorem on page 753 in  \cite{S}
that  $\epsilon^{(n)}(\tau)$ is generic.
The local  uniqueness of the Whittaker function 
is proved in \cite{S}. 

There is one example where one can use the theory of small representations to study the above Conjecture. This case is for the group $SL_4^{(2)}({\bf A})$ or the group $SO_6^{(2)}({\bf A})$. Since these groups are closely related to $GL_4^{(2)}({\bf A})$, one can hope that proving a similar conjecture for the first two, would imply the conjecture for 
$GL_4^{(2)}({\bf A})$.

For the group $SL_4^{(2)}({\bf A})$, one can use the fact that $SL_4$ and $SL_2$ are a commuting pair in $F_4$. Then use the minimal representation on $F_4^{(2)}({\bf A})$ to establish a similar Conjecture for $SL_4$. The global setup for this is explained in details in \cite{G2} Section 3.3. In that reference, the local correspondence is not established. The second construction involves the dual pair $SO_3\times SO_6$ in the group $SO_9$. Using the representation constructed in \cite{B-F-G} on the group $SO_9^{(2)}({\bf A})$, one can establish a similar Conjecture for this case. We hope to study this case in the near future, as  related to a more general type of construction on metaplectic covering groups.

\subsection{The Global Construction 1: The Whittaker Integral}\label{suz2}
Let $\pi^{(n)}$ denote an irreducible cuspidal representation defined on the group $GL_r^{(n)}({\bf A})$. Let $\tau$ denote the cuspidal representation of $GL_m({\bf A})$ as was defined in Section \ref{suz1}.
We  assume that $r\le m$.
In this Section we establish how, using Conjecture \ref{conj1}, one can construct global integrals which unfold to the local integrals \eqref{ra2}. We assume the existence of a representation $\epsilon^{(n)}(\tau)$ as in Conjecture \ref{conj1}. 

Let $Q_{r+1}$ denote the parabolic subgroup of $GL_{nm}$ whose Levi part is given by $GL_{r+1}\times GL_1\times GL_1\times\ldots\times GL_1$ where $GL_1$ occurs $nm-r-1$. Since $r\le m$ and $n>1$, then clearly $nm-r-1\ge 1$. Let $U_{r+1}$ denote the unipotent radical of $Q_{r+1}$. We define a character $\psi_{U_{r+1}}$ on the group $U_{r+1}(F)\backslash U_{r+1}({\bf A})$ as follows. For $u=(u_{i,j})\in U_{r+1}$ define $\psi_{U_{r+1}}(u)=\psi(u_{r+1,r+2}+u_{r+2,r+3}+\cdots +u_{nm-1,nm})$.

There are two cases to consider. 
Assume first, that either $r<m$ or that if $r=m$, then the representation 
$\epsilon^{(n)}(\tau)$ is cuspidal. In this case we define the following global integral.
\begin{equation}\label{glob4}
\int\limits_{GL_r(F)\backslash GL_r({\bf A})}
\int\limits_{U_{r+1}(F)\backslash U_{r+1}({\bf A})}
\varphi_{\pi^{(n)}}(g)\overline{\varphi}_{\epsilon^{(n)}(\tau)}(ug_0)
\psi_{U_{r+1}}(u)|g|^{s-\frac{nm-r}{2}}dudg
\end{equation}
Here, for $g\in GL_r$, we denote $g_0=\text{diag}(g,I_{nm-r})$. This integral converges absolutely for all $s$, even in the case when 
$\epsilon^{(n)}(\tau)$ is the residue representation. This follows from the fact that on the center of $GL_r$ the residue representation 
$\epsilon^{(n)}(\tau)$ is rapidly decreasing. Indeed, this last statement holds  since the constant term of 
$\epsilon^{(n)}(\tau)$ along the unipotent radical $L_r$ is zero. Here $L_r$ is the unipotent radical of the maximal parabolic subgroup of $GL_{nm}$ whose Levi part is $GL_r\times GL_{nm-r}$. It is clear that the unfolding process of integral \eqref{glob4} is the same as the unfolding process of the well known Rakin-Selberg integral ( see \cite{J-PS-S} ) which represents the $L$ function for $GL_r\times GL_{nm}$. From this we obtain that for $\text{Re}(s)$ large integral \eqref{glob4} is equal to
\begin{equation}\label{glob5}
\int\limits_{ V_r({\bf A})\backslash GL_r({\bf A})}W_{\pi^{(n)}}(g)\overline{W}_{\epsilon^{(n)}(\tau)}\begin{pmatrix} g&\\
&I_{nm-r}\end{pmatrix}|g|^{s-\frac{nm-r}{2}}dg
\end{equation}
From this, and from the local computations of the proceeding sections, we obtain
\begin{theorem}\label{th2}
Assume that Conjecture \ref{conj1} holds. Assume that $r<m$ or that if $r=m$ then $\tau$ is not in the image of the Shimura lift from $GL_m^{(n)}({\bf A})$. Let $S$ denote a finite set of places including the archimedean places, such that outside $S$ all data is unramified. Then the partial $L$ function $L^S(\pi^{(n)}\times \tau,s)$ is a holomorphic function for all $s$.
\end{theorem}
Notice that if $\tau$ is  in the image of the Shimura lift of $\tau^{(n)}$, then 
$L^S(\pi^{(n)}\times \tau,s)=L^S(\pi^{(n)}\times \tau^{(n)},s)$.

The second case to consider is when $r=m$ and that $\epsilon^{(n)}(\tau)$ is the residue representation. As was mentioned above the existence of $\epsilon^{(n)}(\tau)$ is not known. Henceforth, we will assume that it exists.
Since $r=m$, this last representation is defined on the group $GL_{nr}^{(n)}({\bf A})$.

We now construct the global integral in a similar way as in \cite{B-F} integral (2.1). Let
\begin{equation}\label{glob6}
L_{\varphi_{\epsilon^{(n)}(\tau)^{(n)}}}(g_0)=\int\limits_{U_{r+1}(F)\backslash U_{r+1}({\bf A})}
\overline{\varphi}_{\epsilon^{(n)}(\tau^{(n)})}(ug_0)
\psi_{U_{r+1}}(u)du-\int\limits_{U_{r}(F)\backslash U_{r}({\bf A})}
\overline{\varphi}_{\epsilon^{(n)}(\tau^{(n)})}(ug_0)
\psi_{U_{r},0}(u)du
\end{equation}
Here, the character $\psi_{U_{r},0}$ is the trivial  extension of 
$\psi_{U_{r+1}}$ from the the group $U_{r+1}$ to $U_r$.
Then we consider the integral
\begin{equation}\label{glob7}
I=\int\limits_{GL_r(F)\backslash GL_r({\bf A})}
\varphi_{\pi^{(n)}}(g)
L_{\varphi_{\epsilon^{(n)}(\tau^{(n)})}}(g_0)
|g|^{s-\frac{nm-r}{2}}dg
\end{equation}
where we recall that $g_0=\text{diag}(g,I_{r(n-1)})$. 

With this construction we can now prove,
\begin{theorem}\label{th3}
Let $\pi^{(n)}$ and  $\tau^{(n)}$ denote two irreducible cuspidal representations  of the group $GL_r^{(n)}({\bf A})$. Assume that the residue representation $\epsilon^{(n)}(\tau^{(n)})$ exist. Then the partial $L$ function $L^S(\pi^{(n)}\times \tau^{(n)},s)$ has a meromorphic continuation to the whole complex plane. If $\pi^{(n)}$ is not isomorphic to the  representation $\widehat{\tau}^{(n)}$, then this partial $L$ function is holomorphic. If $\pi^{(n)}$ is  isomorphic to $\widehat{\tau}^{(n)}$ then the partial $L$ function can have at most a simple pole at $s=0$ and $s=1$.
\end{theorem}
The representation $\widehat{\tau}^{(n)}$ was defined right after Theorem \ref{mainth}.
We remark that one expects a stronger result. Indeed, it is conjectured that the partial $L$ function $L^S(\pi^{(n)}\times \tau^{(n)},s)$ will have a pole at $s=1$ if and only if $\pi^{(n)}$ is isomorphic to the representation
$\widehat{\tau}^{(n)}$.
\begin{proof}
First, for $\text{Re}(s)$ large we perform the same Fourier expansions as in the non-covering case, see \cite{J-PS-S}, and we obtain that the integral $I$ is equal to integral \eqref{glob5} with $\epsilon^{(n)}(\tau^{(n)})$ instead of $\epsilon^{(n)}(\tau)$. Since data can be chosen that at the bad places the corresponding local integral is no zero, it follows from Theorem \ref{th1}, that there is a choice of data such that up to a factor which is holomorphic in $s$, integral $I$ is equal to the partial $L$ function $L^S(\pi^{(n)}\times\tau^{(n)},ns-(n-1)/2)$.

Next we consider the continuation of integral $I$. This is done as in the case when $n=1$. See \cite{J-PS-S}.
Starting with integral \eqref{glob7}, we factor the integration domain to $$I=\int\limits_0^\infty\ \int\limits_{GL_r(F)\backslash GL_r^0({\bf A})}
=\int\limits_0^1 \int\limits_{GL_r(F)\backslash GL_r^0({\bf A})}+
\int\limits_1^\infty\ \int\limits_{GL_r(F)\backslash GL_r^0({\bf A})}$$
where $GL_r^0({\bf A})=\{g\in GL_r({\bf A})\ :\ |\text{det} g|_{\bf A}=1\}$. Denote the right hand side as $I_1+I_2$, where $I_2$ is the right most summand. It follows from the fact that $L_{\varphi_{\epsilon^{(n)}(\tau)^{(n)}}}(t_0g_0)$ is rapidly decreasing in $t$ in that domain, that integral $I_2$ converges absolutely for all $s$. Here $t_0=\text{diag}(tI_r,I_{r(n-1)})$. 

In integral $I_1$ change variables $t\to t^{-1}$. Thus, integral $I_1$ is equal to
\begin{equation}\label{glob8}
I_1=\int\limits_1^\infty\int\limits_{GL_r(F)\backslash GL_r^0({\bf A})}
\varphi_{\pi^{(n)}}(g)
L_{\varphi_{\epsilon^{(n)}(\tau^{(n)})}}(t_0^{-1}g_0)
\omega_{\pi^{(n)}}(t^{-1})t^{-rs+\frac{r^2(n-1)}{2}}d^*tdg
\end{equation}
Here $\omega_{\pi^{(n)}}$ is the central character of $\pi^{(n)}$. From equation \eqref{glob6} it follows that $L_{\varphi_{\epsilon^{(n)}(\tau^{(n)})}}(t_0^{-1}g_0)$ is a sum of two terms. The first term is equal to
\begin{equation}\label{glob9}
\int\limits_{U_{r+1}(F)\backslash U_{r+1}({\bf A})}
\overline{\varphi}_{\epsilon^{(n)}(\tau^{(n)})}(ut_0^{-1}g_0)
\psi_{U_{r+1}}(u)du=
\int\limits_{U_{r+1}(F)\backslash U_{r+1}({\bf A})}
{\varphi}_{\epsilon^{(n)}(\widehat{\tau}^{(n)})}(^tut_0(^tg_0^{-1}))\psi_{U_{r+1}}(u)du
\end{equation}
Here, $^th$ denotes the transpose matrix of $h$. We used the fact that the space of functions ${\varphi}_{\epsilon^{(n)}(\widehat{\tau}^{(n)})}$  can be realized as the space of functions $\varphi_{\epsilon^{(n)}(\tau^{(n)})}(^th^{-1})$. 

Let $Y$ denote the unipotent subgroup of $GL_{nr}$ consisting of all matrices of the form
$$y=\begin{pmatrix} I_r&&\\ y'&I_{rn-r-1}&\\ &&1\end{pmatrix}\ \ 
\ \ \ \ \ y'\in Mat_{(rn-r-1)\times r}$$
Let $X_{r,n}$ denote the maximal unipotent subgroup of $GL_{rn-r-1}$ consisting of upper unipotent matrices. We embed $X_{r,n}$ inside $GL_{nr}$ as all matrices of the form $\text{diag} (I_{r+1},x)$ with $x\in X_{r,n}$. Denote $w_r=\text{diag}(I_{r+1},J_{rn-r-1})$.

Then, the right hand side of equation \eqref{glob9} is equal to
\begin{equation}\label{glob10}
\int\limits_{Y(F)\backslash Y({\bf A})}
\int\limits_{X_{r,n}(F)\backslash X_{r,n}({\bf A})}
(w_r{\varphi})_{\epsilon^{(n)}(\widehat{\tau}^{(n)})}(xyt_0(^tg_0^{-1}))\psi_{X_{r,n}}(x)dxdy
\end{equation}
Here $\psi_{X_{r,n}}$ is defined as the restriction of $\psi_{U_{r+1}}$ to the group $X_{r,n}$. As in \cite{J-PS-S}, it follows from simple application of root exchange, see \cite{G-R-S} for this process, that integral \eqref{glob10} is equal to
\begin{equation}\label{glob11}
\int\limits_{Y({\bf A})}
\int\limits_{U_{r+1}(F)\backslash U_{r+1}({\bf A})}
(w_r{\varphi})_{\epsilon^{(n)}(\widehat{\tau}^{(n)})}(uyt_0(^tg_0^{-1}))\psi_{U_{r+1}}(u)dudy
\end{equation}
Thus, integral $I_1$ is equal to
\begin{align}\label{glob12}
\int \varphi_{\pi^{(n)}}(g) [
\int\limits_{Y({\bf A})}
\int\limits_{U_{r+1}(F)\backslash U_{r+1}({\bf A})}
(w_r{\varphi})_{\epsilon^{(n)}(\widehat{\tau}^{(n)})}(uyt_0(^tg_0^{-1}))\psi_{U_{r+1}}(u)dudy -\\
- \int\limits_{U_{r}(F)\backslash U_{r}({\bf A})}
\overline{\varphi}_{\epsilon^{(n)}(\tau^{(n)})}(ut_0^{-1}g_0)
\psi_{U_{r},0}(u)du ]\ \omega_{\pi^{(n)}}(t^{-1})t^{-rs+\frac{r^2(n-1)}{2}}d^*tdg
\end{align}
Here the outer integration domain is as in integral \eqref{glob8}. Inside the parentheses  add and subtract the term
\begin{equation}\label{glob13}
\int\limits_{Y({\bf A})}
\int\limits_{U_{r}(F)\backslash U_{r}({\bf A})}
(w_r{\varphi})_{\epsilon^{(n)}(\widehat{\tau}^{(n)})}(uyt_0(^tg_0^{-1}))\psi_{U_{r},0}(u)dudy
\end{equation}
Then integral \eqref{glob12} is a sum of three integrals. The first is 
\begin{equation}\label{glob14}\notag
\int_1^\infty\int \varphi_{\pi^{(n)}}(g)[ 
\int\limits_{Y({\bf A})}
\int\limits_{U_{r+1}(F)\backslash U_{r+1}({\bf A})}
(w_r{\varphi})_{\epsilon^{(n)}(\widehat{\tau}^{(n)})}(uyt_0(^tg_0^{-1}))\psi_{U_{r+1}}(u)dudy -
\end{equation}
$$-\int\limits_{U_{r}(F)\backslash U_{r}({\bf A})}
(w_r{\varphi})_{\epsilon^{(n)}(\widehat{\tau}^{(n)})}(ut_0(^tg_0^{-1}))\psi_{U_{r},0}(u)du ]\ \omega_{\pi^{(n)}}(t^{-1})t^{-rs+\frac{r^2(n-1)}{2}}d^*tdg$$
As in the case of integral $I_2$, this integral converges absolutely for all $s$. The second integral to consider is 
\begin{equation}\label{glob15}
\int_1^\infty\int 
\int\limits_{Y({\bf A})}
\int\limits_{U_{r}(F)\backslash U_{r}({\bf A})}
\varphi_{\pi^{(n)}}(g)
(w_r{\varphi})_{\epsilon^{(n)}(\widehat{\tau}^{(n)})}(uyt_0(^tg_0^{-1}))\psi_{U_{r},0}(u)\ \omega_{\pi^{(n)}}(t^{-1})
t^{-rs+\frac{r^2(n-1)}{2}}dudy d^*tdg
\end{equation}
Conjugate $t_0(^tg_0^{-1})$ to the left. First we obtain a factor of $t^{-r(r(n-1)-1)}$ from the change of variables in $Y$. Then we need to apply the constant term which is a part of the integration over $U_r$. This gives us a factor of 
\begin{equation}\label{glob16}
\omega_{\widehat{\tau}^{(n)}}(t)|tI_r|^{-\frac{n-1}{2n}}\delta_{P_{n,r}}^{1/2}\begin{pmatrix} tI_r&\\
&I_{r(n-1)}\end{pmatrix}=\omega_{\widehat{\tau}^{(n)}}(t)t^{-\frac{r(n-1)}{2n}+\frac{r^2(n-1)}{2}}
\end{equation}
Here $\omega_{\widehat{\tau}^{(n)}}$ is the central character of $\widehat{\tau}^{(n)}$.
In the $^tg_0^{-1}$ variable, applying the constant term,  we argue in a similar way as in \cite{J-L} to deduce that we get a certain vector in the space of $\widehat{\tau}^{(n)}$ evaluated at $^tg^{-1}$. Plugging all this into integral \eqref{glob15}, we obtain
\begin{equation}\label{glob17}
\int\limits_{GL_r(F)\backslash GL_r^0({\bf A})}
\varphi_{\pi^{(n)}}(g){\varphi}_{\widehat{\tau}^{(n)},1}(^tg^{-1})dg
\int\limits_1^\infty \omega_{\pi^{(n)}}^{-1}\omega_{\widehat{\tau}^{(n)}}(t)t^{-r(s-1+\frac{n-1}{2n})}d^*t
\end{equation}
Here $\varphi_{\widehat{\tau}^{(n)},1}$ is a certain vector in the space of $\widehat{\tau}^{(n)}$. The product in equation \eqref{glob17} is a meromorphic function in $s$, and has a simple pole at $s_0=\frac{n+1}{2n}$ if and only if $\pi^{(n)}=\widehat{\tau}^{(n)}$. 
Notice that $ns_0-\frac{n-1}{2}=1$. This implies that the partial $L$ function $L^S(\pi^{(n)}\times\tau^{(n)},ns-(n-1)/2)$ can have a simple pole at $s=\frac{n+1}{2n}$.

Finally, the last integral to consider is
\begin{equation}\label{glob18}
\int_1^\infty\int 
\int\limits_{U_{r}(F)\backslash U_{r}({\bf A})}
\varphi_{\pi^{(n)}}(g)
\overline{\varphi}_{\epsilon^{(n)}(\tau^{(n)})}(ut_0^{-1}g_0)
\psi_{U_{r},0}(u)\ \omega_{\pi^{(n)}}(t^{-1})t^{-rs+\frac{r^2(n-1)}{2}}dud^*tdg
\end{equation}
Conjugating $t_0^{-1}$ to the left, we obtain in a similar way as in equation \eqref{glob16}, the factor 
\begin{equation}\label{glob19}\notag
|t^{-1}I_r|^{-\frac{n-1}{2n}}\delta_{P_{n,r}}^{1/2}\begin{pmatrix} t^{-1}I_r&\\ &I_{r(n-1)}\end{pmatrix}=t^{\frac{r(n-1)}{2n}-\frac{r^2(n-1)}{2}}
\end{equation}
Arguing as above, we obtain that integral \eqref{glob18} can have at most a simple pole at $s_0=\frac{n-1}{2n}$, and this can happen only if $\pi^{(n)}=\widehat{\tau}^{(n)}$.

\end{proof}

\subsection{The Global Construction 2: The Doubling Integral }\label{suz20}
In this section we introduce a global doubling integral. We do it in the case when $n=1$. It is possible that a certain variation of this construction can be applied to the metaplectic covering groups. The connection of this doubling integral to the content of the previous sections of this paper, is the fact that it unfolds to an Eulerian integral whose local factors are the integrals \eqref{ra3}. From this we deduce that this integral represents the tensor product partial $L$ function $L^S(\pi\times\tau,s)$.

In details, let $\pi$ denote an irreducible cuspidal representation of $GL_r({\bf A})$, and let $\tau$ denote an irreducible cuspidal representation of $GL_m({\bf A})$. Let
${\mathcal E}_\tau$ denote the Speh representation defined on the group $GL_{rm}({\bf A})$, which is attached to the cuspidal representation $\tau$. When $n=1$, this is the residue representation $\epsilon^{(n)}(\tau)$ defined in sub-section \ref{suz1}.
In contrast to the cases when $n>1$,  the existence of ${\mathcal E}_\tau$ is well known. 

In \eqref{mat1} we introduced the unipotent group $U_{nm,r}$. We denote by $U_{m,r}$ the group $U_{nm,r}$ with $n=1$.
Let  $U'_{m,r}$ denote the subgroup of $U_{m,r}$ consisting of all matrices $u\in U_{m,r}$ as in equation \eqref{mat1} with the condition that $X_{1,2}=0$. The character $\psi_{U_{m,r}}$ defined right after equation \eqref{mat1} defines by restriction a character of $U'_{m,r}$. We denote this character by $\psi_{U_{m,r}'}$. Consider the Fourier coefficient 
\begin{equation}\label{glob21}
f_1(h)=\int\limits_{U'_{m,r}(F)\backslash U'_{m,r}({\bf A})}E_\tau(u'h)\psi_{U_{m,r}'}(u')du'
\end{equation}
Here, $E_\tau$ is a vector in the space of ${\mathcal E}_\tau$, and $h\in GL_{rm}({\bf A})$. 

We will also need the function
\begin{equation}\label{glob22}
f_2(h)=\int\limits_{U_{m,r}(F)\backslash U_{m,r}({\bf A})}E_\tau(uh)\psi_{U_{m,r}'}(u)du
\end{equation}
Here, the character $\psi_{U_{m,r}'}$ is the trivial extension from the group $U'_{m,r}$ to $U_{m,r}$. Notice that integral \eqref{glob22} contains as an inner integration, the constant term along the unipotent radical of the maximal parabolic subgroup of 
$GL_{rm}$ whose Levi part is $GL_r\times GL_{r(m-1)}$.

Embed the group $GL_r\times GL_r$ inside $GL_{rm}$ as 
\begin{equation}\label{glob23}
(g_1,g_2)\to \text{diag}(g_1,g_2,\ldots,g_2)
\end{equation}
where $g_2$ appears $m-1$ times. Let $Z^\Delta$ denote the center of
$GL_r$ embedded diagonally in $GL_r\times GL_r$. Then $Z^\Delta$ is also the center of $GL_{rm}$.

The global integral we consider is
\begin{equation}\label{glob24}
\int\limits_{Z^\Delta({\bf A})(GL_r(F)\times GL_r(F))\backslash (GL_r({\bf A})\times GL_r({\bf A}))}\varphi_{1}(g_1)\varphi_{2}(g_2)[f_1((g_1,g_2))-f_2((g_1,g_2))]\left |\frac{g_1}{g_2}\right |^{s'}dg_1dg_2
\end{equation}
Here $s'$ is defined right after equation \eqref{ra3}. Also, the functions $\varphi_1$ and $\varphi_2$ are two vectors in the space of $\pi$. 
As in integral \eqref{glob7}, this integral converges for $Re(s)$ large, and admits a meromorphic continuation to the whole complex plane.
We now show that this integral unfolds to an integral whose local components are given by integral \eqref{ra3}. 

Starting with the integral $f_1((g_1,g_2))$, we consider its Fourier expansion along the subgroup $U_0$ of $U_{m,r}$ defined by all matrices as in equation \eqref{mat1} such that $X_{i,j}=0$ for all $(i,j)\ne (1,2)$. Thus $U_{m,r}=U_0U'_{m,r}$. The group $GL_r(F)\times GL_r(F)$ acts on this expansion 
with a finite number of orbits. 

The contribution to integral \eqref{glob24} from the trivial character in the above expansion is zero since it is canceled with $f_2((g_1,g_2))$. The contribution from all other orbits except the full rank orbit, is zero.
This follows from the cuspidality of $\pi$. Indeed, in this case the stabilizer inside all matrices $g_1\in GL_r$ contains a unipotent radical of a maximal parabolic subgroup of $GL_r$. We denote this unipotent radical by $X$. Further Fourier expansion, it follows from Proposition 2.4 in \cite{G-S} that we obtain as inner integration the integral  $\int\varphi_1(xg_1)dx$. Here the integration domain is $X(F)\backslash X({\bf A})$. Thus we get zero contribution by the cuspidality of $\pi$. 

We are left with the full rank orbit whose stabilizer inside $GL_r\times GL_r$ is the diagonal embedding of $GL_r$ in $GL_r\times GL_r$. By the properties of the Speh representation, the corresponding Fourier coefficient is left invariant under the diagonal embedding of the group $GL_r({\bf A})$ inside $GL_r({\bf A})\times GL_r({\bf A})$.  Thus, integral \eqref{glob24} is equal to 
\begin{equation}\label{global25}
\int\limits_{ GL_r}\omega_{\pi}(g)W_{\tau,rm}\begin{pmatrix} g&\\
&I_{r(m-1)}\end{pmatrix}|g|^{s'}dg
\end{equation}
Here, 
\begin{equation}\label{global26}\notag
\omega_{\pi}(g)=\int\limits_{Z({\bf A})GL_r(F)\backslash GL_r({\bf A})}
\varphi_1(h)\varphi_2(hg)dh
\end{equation}
and 
\begin{equation}\label{global27}\notag
W_{\tau,rm}(h')=\int\limits_{U_{m,r}(F)\backslash U_{m,r}({\bf A})}
E_\tau(uh')\psi_{U_{m,r}}(u)du
\end{equation}
From this it follows that integral \eqref{global25} is factorizable, and its corresponding local integrals are given by integrals \eqref{ra3}.

\end{document}